\documentclass[12pt,twoside]{article}

\usepackage{amsfonts}
\usepackage{bezier}
\usepackage{euscript}
\usepackage{amsfonts,amssymb,amsmath,amsbsy,amsthm,amscd}
\usepackage[T2A]{fontenc}
\usepackage[cp1251]{inputenc}
\usepackage[english,russian]{babel}   

\usepackage[russian]{babel}
\usepackage{enumitem}

\newtheorem{theorem}{Теорема}

\begin{document}

\def\udk{519.766}
\title{КРИТЕРИЙ ПОДСТАНОВОЧНОСТИ ПАЛИНДРОМОВ ШТУРМА}
{И.\,А.~Решетников,\footnote[1]{%
{\it Решетников Иван Андреевич} --- аспирант ФИВТ МФТИ(ГУ),
e-mail: reshetnikov.ivan@phystech.edu.}
А.\,Я.~Канель-Белов\footnote[2]{%
{\it Алексей Яковлевич Канель-Белов} --- Shengeng University, BIU
e-mail: kanelster@gmail.com.}
}

\begin{abstract}
  В статье приводится критерий подстановочности симметричных бесконечных в обе стороны слов Штурма и его доказательство.

The article provides a criterion for the substitution of symmetric Sturm words infinite on both sides and its proof. This work was carried out with the help of the Russian Science Foundation Grant N 17-11-01377.

\end{abstract}

{\bf keywords:}
Sturm's words, substitutive words, mechanical words, dynamical systems, circle rotation, Rausy induction, symbolic dynamics

\section{Введение}

\emph{Двусторонне-бесконечное слово} $w$ --- это отображение $\mathbb{Z} \to A$,
где $A$ --- алфавит слова $w$. Будем обозначать $w(n)$ через $w_n$.

Пусть есть символическая динамика $(M, R_\alpha, x_0, U)$, где $M$ --- окружность,
$U$ --- дуга угловой меры $\alpha$, ($\alpha$ иррациональное),
$R_\alpha$ --- поворот на $\alpha$ относительно центра окружности
(угловую меру всей окружности считаем единицей) ---  функция эволющии,
 $x_0$ --- начальная точка --- начало дуги $\alpha$.
 Динамическая система $(M, R_\alpha)$ порождает некоторое бесконечное слово $w$,
  эволюцию точки $x_0$. Такие слова называются механическими.

Двусторонне-бесконечное слово $w$ называется \emph{палиндромом}, если
$\forall n \in \mathbb{Z}: \; w_n=w_{-n}$
(палиндром, симметричный относительно буквы) или
$\forall n \in \mathbb{Z}: \; w_n=w_{1-n}$ (палиндром, симметричный относительно междубуквия).

Правосторонне-бесконечное слово $w$ над алфавитом $A$ называется \emph{чисто подстановочным},
если оно представляется в виде $w=\phi^\infty(a)$, где $a \in A$ --- буква,
а $\phi$ --- подстановка, такая, что $\phi(a)=aU$, где $U \in A^\ast$, непустое.

Слово $w$ называется подстановочным, если получается из чисто подстановочного
слова $w'$ подстановкой $h$ применённой к слову $w'$: $w=h(w')$.

Тогда $\phi( \phi(a) )=aU\phi(U)$, $\phi( \phi( \phi(a) ) )=aU\phi((U\phi(U))$
и т.д. Предельным переходом получаем $w=\phi^\infty(a)$.

Понятия \emph{подстановочности и чистой подстановочности слов, бесконечных
в обе стороны} вводятся аналогично с той лишь разницей, что для двусторонне
бесконечного подстановочного слова подстановка должна
обладать свойством $\phi(a)=Ua$ и $\phi(b)=bW$. Тогда двусторонне-бесконечное
слово $w$
 будет конкатенацией левосторонне-бесконечного подстановочного слова
 (позиции от $-\infty$ до 0 включительно)и правосторонне-бесконечного
 подстановочного слова (позиции большие нуля) и будет представимо в виде $w=\phi^\infty(a|b)$.

 Мотивировки задач и подробности --- см. [1], [2].

\section{Индукция Рози}

Метод индукции Рози [3] позволяет перейти от изучения некоторых
механических слов к изучению цепных дробей. Его можно сформулировать
следующим образом.

Пусть есть символическая динамика $(M, R_\alpha, x_0, U)$, где $M$ --- окружность,
$U$ --- дуга угловой меры $\alpha$, ($\alpha$ иррациональное),
$R_\alpha$ --- поворот на $\alpha$ относительно центра окружности
(угловую меру всей окружности считаем единицей) ---  функция эволющии,
$x_0$ --- начальная точка --- начало дуги $\alpha$.
 Динамическая система $(M, R_\alpha)$ порождает некоторое бесконечное
 слово $w$, эволюцию точки $x_0$.

Требование иррациональности числа $\alpha$ продиктовано периодичностью
сдвигов на рациональный угол, этот тривиальный случай мы не разбираем.

Метод индукции Рози состоит в следующем. Мы преобразуем символическую
динамику $(M, R_\alpha, x_0, U)$ в символическую динамику
$(\tilde{M}, \tilde{R_\alpha}, \tilde{x_0}, \tilde{U})$ по следующим правилам.

\begin{enumerate}

\item
Если $\alpha < \dfrac{1}{2}$, то в слове $w$ после $a$ всегда следует $b$.
Поэтому слово $w$ получается из слова $\tilde{w}$ заменой
$
   a \to ab;
   b \to b.
$
Возьмём $\tilde{M}$ --- окружность длины $1-\alpha$. Более наглядно, мы
вырезаем из окружности $M$ дугу длины $\alpha$ (следующую за выделенной
дугой $U$). Тогда слово $\tilde{w}$ будет порождаться символической динамикой
$(\tilde{M}, \tilde{R_\alpha}=R_\alpha, \tilde{x_0}=x_0, U=\tilde{U})$.

\item
Если же $\alpha > \dfrac{1}{2}$, то заменим $a$ на $b$ и наоборот.
Для этого надо положить $\tilde{U}=M\backslash U$.

\end{enumerate}

Можно заметить, что описанный алгоритм аналогичен алгоритму разложения
 $\alpha$ в цепную дробь. Таким образом, если $\alpha$ --- квадратичная
 иррациональность (её цепная дробь периодична), то в какой-то мы получим
 символическую динамику, эквивалентную уже встречавшейся, а значит можем
  записать подстановку, с помощью которой можно вырастить слово $w$.
  Обозначим за $\varphi$ композицию подстановок, образующих период,
  а за $\psi$ комозицию подстановок, образующих предпериод. Слово $w$
  представимо в виде $w=\psi \circ \varphi^\infty(a)$, то есть подстановочно.
  Если же предпериод отсутствует, то слово $w$ чисто подстановочно.

Если же $\alpha$ --- не квадратичная иррациональность, то процесс не
зациклится, но в любом случае можно будет изучать свойства цепной дроби
 $\alpha$ и соотносить их со свойствами слова $w$.

\section{Критерий подстановочности палиндромных слов Штурма}

Всего есть три двусторонне бесконечных палиндрома, являющихся механическими
словами параметра $\alpha$. Это слова, отвечающие серединам дуг и слово,
получающееся из точки, которая первым поворотом переходит в точку,
симметричную ей относительно оси симметрии динамической системы. Палиндромы,
получающиеся из середин дуг симметричны относительно своей буквы, а третий
палиндром симметричен относительно междубуквия.

\subsection{Аналог индукции Рози для палиндромов}

В этом разделе мы докажем следующий критерий подстановочности палиндромов
Штурма. Во всём разделе полагаем алфавит $A=\{0, 1\}$.

\begin{theorem}[Критерий подстановочности палиндромов Штурма]
Символическая динамика $(M, R_\alpha, x_m, U)$, где $x_m$ --- середина дуги
$U$ длины $\alpha$ на окружности $M$ длины 1 порождает подстановочный
двусторонне-бесконечный палиндром $p$ тогда и только тогда, когда $\alpha$
--- квадратичная иррациональность
\end{theorem}

Так как наше преобразование окружности обратимо, то можно говорить и о
двусторонне-бесконечных словах, порождаемых символической динамикой.
В частности, о словах-палиндромах. Доказательство будет основано на
сходстве некоторого алгоритма с индукцией Рози.

{\bf Доказательство.}
В одну сторону: если слово $p$ является подстановочным, то $\alpha$
--- квадратичная иррациональность. $\alpha$ -- доля числа единиц.
Пусть подстановка $\phi$ порождает слово $p$. Пусть при действии порождающей
 подстановки $\phi$ 0 переходит в слово с $a$ нулями и $b$ единицами,
 а 1 переходит в слово с $c$ нулями и $d$ единицами. Тогда,
 так как $p$ --- неподвижная точка подстановки $\phi$, получаем
 квадратичное относительно $\alpha$ уравнение
$$
\alpha = \dfrac{b(1-\alpha)+d \alpha}{(a+b)(1-\alpha)+(c+d)\alpha},
$$
значит, $\alpha$ --- квадратичная иррациональность.

В другую сторону, если $\alpha$ --- квадратичная иррациональность,
то $p$ подстановочно. Пусть есть бесконечное слово Штурма $p$ ---
палиндром над алфавитом $\{0, 1\}$. Пусть в нём доля единиц равна $\alpha$.
Если нулей меньше, чем единиц, то сделаем подстановку $E:0\leftrightarrow1$.
 Остаётся расписать действия алгоритма при $\alpha<\dfrac{1}{2}$.
 Так как $p$ является словом Штурма и $\alpha<\dfrac{1}{2}$ в нём
 не могут встретиться две единицы подряд. Возможны несколько случаев.

\begin{enumerate}
\item Слово $p$ симметрично относительно 1. Будем относить эту 1 к
левой половине слова $p$. Так как в слове $p$ нулей больше, чем
единиц, то перед каждой единицей идёт ноль. Значит, слово $p$ получается
 подстановкой
$
G:    1 \to 01;
   0 \to 0\\
$
из некоторого слова $p'$. Причём подстановка подобрана так, что левая
и правая части слова $p$ получаются из левой и правой части слова $p'$.
Так как количество нулей в каждой группе нулей уменьшилось на 1, а больше
 ничего не произошло, то слово $p'$ также является палиндромом,
 симметричным относительно 1.

\item
Слово $p$ симметрично относительно 0. Будем относить этот 0 к левой
 половине слова~$p$. Так как в слове $p$ нулей больше, чем единиц,
 то после каждой единицы идёт ноль. Значит, слово $p$ получается
 подстановкой
$
\tilde{G}:
   1 \to 10;
   0 \to 0
$
из некоторого слова $p'$. Левая и правая части слова $p$ получаются
соответсвенно из левой и правой частей слова $p'$. Количество нулей
 в каждой группе нулей уменьшилось на 1, поэтому слово $p'$ также
 является палиндромом, но на этот раз симметричным относительно
  междубуквия.

\item
Слово $p$ симметрично относительно междубуквия. Тогда буквы $w_0$
и $w_1$ должны быть нулями. Тогда слово $p$ получается подстановкой
$G$ из некоторого слова $p'$, при этом слово $p'$ симметрично
относительно нуля.
\end{enumerate}

В нашем случае $\alpha$ --- квадратичная иррациональность.
При действии приведённого алгоритма $\alpha$ претерпевает такие
 же изменения, как и в индукции Рози, а значит, так как $\alpha$
 --- квадратичная иррациональность, то доля единиц $\alpha$ начиная
 с некоторого момента будет изменяться периодически. Для каждого
 $\alpha$ возможны лишь три варианта симметричности, три различных
  палиндрома. Поэтому возможно лишь конечное число пар $(\alpha, j)$,
  где $j$ --- номер варианта симметричности. Для каждой пары
  однозначно определён переход согласно алгоритму к некоторой другой паре.
  Поэтому процесс и в общем также цикличен. Из цикличности
  процесса следует, что палиндром $p$ подстановочен: как и в случае с обычной индукцией Рози.
Теорема доказана.

\section{Чистая подстановочность двусторонне-бесконечного палиндрома Фи\-бо\-нач\-чи-Штур\-ма}

Если мы предъявим подстановку $\phi$, оставляющую на месте слово $w$, такую,
что $\phi(a)=Ua$ и $\phi(b)=bW$ для каких-то непустых слов $U$ и $W$, тогда
двусторонне-бесконечное слово $w$ будет являться чисто подстановочным.

\begin{theorem}{Утверждение}
Подстановка $\phi : \; 0 \to 00101, 1 \to 001$ оставляет на месте палиндром
 $w$ и обладает свойством $\phi(a)=Ua$ и $\phi(b)=bW$, значит
 двусторонне-бесконечный палиндром $w$ , симметричный относительно 1
 для $\alpha=\phi$ (здесь $\phi$ --- золотое сечение) растится этой
 подстановкой, т.е. $w=\phi^\infty(1|0)$.
\end{theorem}

{\bf Доказательство.}
Данная подстановка $\phi$ получается композицией двух подстановок
$\phi = \psi_2 \circ \psi_1$, где $\psi_1 : \; 0 \to 01, 1 \to 1$
и $\psi_2: \; 0 \to 001, 1 \to 01$, каждая из которых сохраняет отношение
количества нулей и единиц в слове, в котором отношение нулей и единиц
такое же, как в слове фибоначчи. При этом обе эти подстановки сохраняют
свойство палиндромности. Первая подстановка увеличивает количество нулей
 в каждой группе нулей на 1, а значит сохраняет палиндромность. Вторая
 порождает группы нулей, симметричные относительно точки симметрии слова.
  $\psi_1$ слово, симметричное относительно 1, переводит в слово,
  симметричное относительно междубуквия, а $\psi_2$ слово, симметричное
  относительно междубуквия, переводит в слово, симметричное относительно 1.
  Так как такое сбалансированное слово только одно для данного отношения
  количеств нулей и единиц, то $\phi$, переводящая слово, симметричное
  относительно 1 в слово, симметричное относительно 1, при этом оставляя
  неизменным отношение количеств 0 и 1, оставляет слово $w$ без изменений.
Утверждение доказано.

\end{document}